\documentclass[a4paper,12pt]{article}

\usepackage{amsmath}
\usepackage{amsthm}
\usepackage{amssymb}  
\usepackage{latexsym} 
\usepackage[all]{xy}
\usepackage{comment}
\usepackage{rotating}
\usepackage{slashbox}
\usepackage{multirow}
\usepackage{rotating}

\newcommand{\Q}{{\mathbb Q}}

\newcommand{\ord}{\operatorname{ord}}

\newcommand{\tors}{\operatorname{tors}}

\newcommand{\Z}{{\mathbb Z}}

\newfont{\wncyr}{wncyr10 at 12pt}
\newfont{\wncyrten}{wncyr10 at 10pt}

\newenvironment{Proof}{\par\noindent{\sc Proof:}}%
                      {\hspace*{\fill}\nobreak$\Box$\par\medskip}
                       {\hspace*{\fill}\nobreak$\Box$\par\medskip}
\newenvironment{myitemize}
{\begin{itemize}
\setlength{\itemsep}{1pt}
\setlength{\parskip}{0pt}
\setlength{\parsep}{0pt}}
{\end{itemize}}

\newtheorem{Proposition}{Proposition}[section]
\newtheorem{Theorem}[Proposition]{Theorem}
\newtheorem{Lemma}[Proposition]{Lemma}
\newtheorem{Corollary}[Proposition]{Corollary}
\newtheorem{Conjecture}[Proposition]{Conjecture}

\theoremstyle{definition}

\newcounter{nootje}
\begin{document}
\normalsize
\title{On elliptic curves whose conductor is a product of two prime powers}
\author{Mohammad Sadek\\ Department of Mathematics and Actuarial Science\\ American University in Cairo\\ mmsadek@aucegypt.edu}
\date{}
\maketitle
\begin{abstract}{\footnotesize We find all elliptic curves defined over $\Q$ that have a rational point of order $N,\;N\ge 4$, and whose conductor is of the form $p^aq^b$, where $p,q$ are two distinct primes, $a,b$ are two positive integers. In particular, we prove that Szpiro's conjecture holds for these elliptic curves.}
\end{abstract}

\section{Introduction}

Let $E$ be an elliptic curve defined over $\Q$ with minimal discriminant $\Delta_E$. We define the conductor $N_E$ of $E$ to be
\[N_E=\prod_{p\mid \Delta_E}p^{f_p},\;f_p={\ord_p(\Delta_E)-m_E+1},\]
where $m_E$ is the number of components on the special fiber of the N\'{e}ron model of $E$ defined over $\mathbb{F}_p$, $f_p\ge1$, see (\cite{Sil2}, Chapter IV, \S 10, 11). Furthermore, $f_p=1$ if and only if $E$ has multiplicative reduction at $p$. We recall that $N_E$ and $\Delta_E$ have the same prime divisors.

The problem of finding all elliptic curves $E$ defined over $\Q$ of a given conductor has been investigated in many articles. Ogg produced the complete list of elliptic curves whose conductor is a 2-power or $2^a3^b$ , see \cite{Ogg2power} and \cite{Oggsmallconductor}. A series of papers dealt with the problem under the condition that $E$ has a rational torsion point. For example, in \cite{Hadano} the elliptic curves with conductor $p^m$, $p$ is prime, and 2-torsion points were listed.

 It was shown in \cite{Hadanoremarks} that all elliptic curves with conductor $2^m p^n$ where $p\equiv3$ or $5$ mod $8$, $p\ne3$, that have a rational point of order $2$, are effectively determined under the truth of the conjecture of Ankeny-Artin-Chowla.

It is worth mentioning that the complete list of elliptic curves with a prime conductor has already been produced. The following theorem gives this list explicitly.

\begin{Theorem}[Theorem 5.3.2, \cite{Alice}]
\label{thm:primeconductor}
Let $E$ be an elliptic curve over $\Q$ with prime conductor $p$. Then either $|\Delta_E|=p$ or $p^2$, or else $p=11$ and $\Delta_E=11^5$, or $p=17$ and $\Delta_E=17^4$, or $p=19$ and $\Delta_E=19^3$, or $p=37$ and $\Delta_E=37^3$. In particular, $\Delta_E\mid p^5$.
\end{Theorem}
The elliptic curves in Theorem \ref{thm:primeconductor} turn out to satisfy Szpiro's conjecture which is stated below for the convenience of the reader.
\begin{Conjecture}
If $E$ is an elliptic curve over $\Q$, then
\[|\Delta_E|\ll_{\epsilon}N_E^{6+\epsilon}\]
\end{Conjecture}

One of the popular strategies to find elliptic curves $E/\Q$ with a given conductor is to solve certain Diophantine equations obtained by equating the discriminant of $E$ to the product of powers of the prime divisors of the conductor.

Mazur gave a complete classification of the torsion subgroup $E_{\tors}(\Q)$ of $E(\Q)$, see (Theorem 7.5, \S8, Chapter VIII, \cite{sil1}). More precisely, $E_{\tors}(\Q)$ is isomorphic to one of the following fifteen groups:
 \begin{eqnarray*}
 \Z/n\Z,\;1\le n\le 12,\;n\ne11;\textrm{ or }
 \Z/2\Z\times\Z/2n\Z,\;1\le n\le 4.
 \end{eqnarray*}

 Given that $P\in E(\Q)[m],\;m\ne2,3$, it is known that there exist $b,c\in\Q$ such that the following Weierstrass equation defines an elliptic curve $E_{b,c}$ isomorphic to $E$
\begin{eqnarray*}
E_{b,c}:y^2+(1-c)xy-by=x^3-bx^2
\end{eqnarray*}
 with the image of $P$ being $(0,0)$. The discriminant $\Delta(b,c)$ of $E_{b,c}$ is given as follows:
\begin{eqnarray*}
\Delta(b,c)&=&b^3\left(16b^2-b(8c^2+20c-1)-c(1-c)^3\right)
\end{eqnarray*}
By taking $m$ to be an integer in $\{4,5,6,7,8,9,10,12\}$, one finds an explicit relation between $b,c$, see for example \S2 of \cite{LorenziniTamagawa}.

In this article, we generalize Theorem \ref{thm:primeconductor} to elliptic curves whose conductors have two distinct prime divisors only. More precisely, we generate the list of all elliptic curves with $\Q$-rational torsion points of order $N,\;N\ge4,$ whose conductor is a product of two prime powers.

Now we give a brief outline for our approach to solve the problem. As we have seen the family of all elliptic curves with a rational point of order $N$ can be classified using a universal Weierstrass equation. Moreover, we can use certain transformations to write integral Weierstrass equations for the elliptic curves $E_{b,c}$. Since the prime divisors of the minimal discriminant are exactly those of the conductor, we equate the produced discriminant to a product of two prime powers. Consequently, the problem is reduced to solving one or several Diophantine equations.

In fact, the Diophantine equations we produce are more subtle when $N\in\{4,5\}$, whereas the corresponding Diophantine equations are elementary when $N\ge 6$. We collect the harder Diophantine equations in \S \ref{sec:diophantine} for more convenience. Several techniques are followed to attack these equations including elementary methods, factorization over number fields, properties of Lucas sequences, and well-known results from the literature.

Each family of elliptic curves with rational points of order $N$ is treated separately. Given an $N,\,4\le N\le 12,\,N\ne11,$ we list all elliptic curves with a $\Q$-rational $N$-torsion point such that the conductor has only two distinct prime divisors. Moreover, we find a constant $K>0$ such that given such an elliptic curve $E$, the absolute value of the minimal discriminant $\Delta_E$ of $E$ is bounded above by the $K$-th power of the conductor $N_E$. In particular, we prove Szpiro's conjecture for these families of elliptic curves. When $N=10,12$, we show that there are no elliptic curves with an $N$-torsion point whose conductor is a product of two prime powers.

\section{Diophantine equations}
\label{sec:diophantine}

The Catalan's Conjecture (now referred to as Mih\v{a}ilescu's Theorem) will appear several times in this article, so we prefer to state it.
\begin{Proposition}[Mih\v{a}ilescu's Theorem]
\label{prop:catalan}
The only integer solution to the Diophantine equation $\displaystyle x^m-y^n=1$, where $m,n>1$, is $(x,m,y,n)=(\pm3,2,2,3)$.
\end{Proposition}

Now we start solving some Diophantine equations that we will use to prove our main results.

\begin{Lemma}
\label{lem:16pm+1=ql}
There are no integer solutions $(p,m,y,n)$ to the equation
\[16 p^m+1= y^n\]
where $|p|, n$ are prime integers, $m>1$, and $y=l^t$ where $|l|$ is prime and $t>0$.
\end{Lemma}
\begin{Proof}
Let $(p,m,y,n)$ be such solution to $ 16 p^m+1= y^n$. We observe that $y$ is odd. Furthermore, $|p|\ne 2$, otherwise we will have a Catalan's solution $|y^n-2^{m+4}|=1$. So $|p|$ is odd.

We assume $n$ is an odd prime. So $p^my>0$. Since
\[16p^m=y^n-1=(y-1)(y^{n-1}+ y^{n-2}+\ldots+ y+1),\] where the first factor is even, the second factor is the sum of $n$ odd terms, and hence is odd. Therefore $16\mid y-1$.
\begin{myitemize}
\item[i.] If $\displaystyle \gcd(y-1,\frac{y^n-1}{y-1})=1$, then either $y-1=16p^m$ and $\displaystyle \frac{y^n-1}{y-1}=1$, which yields no solutions, or $y=17$ and $\displaystyle\frac{y^n-1}{y-1}=17^{n-1}+\ldots+17+1=p^m$. The latter equation is not solvable for $n\ge 3,m\ge 2$, see Corollary 1 in \cite{BugeaudRoy}. The last possible value $y=-15$ is rejected because it is not a prime power.
\item[ii.] If $p\mid\displaystyle \gcd(y-1,\frac{y^n-1}{y-1})$, then $y\equiv 1$ mod $p$. Moreover, $y-1=\pm 16p^h$, $h>0$, and $\displaystyle \frac{y^n-1}{y-1}= \pm p^{m-h}$. This implies that $n=|p|$ ($n$ is prime). We observe that
    \begin{eqnarray*}
    \pm p^{m-h}&=&\sum_{i=0}^{n-1}(1\pm 16p^h)^i=\sum_{i=0}^{n-1}\sum_{j=0}^i{i\choose j}(\pm16p^h)^j\\
    &=&p+\sum_{i=1}^{n-1}\sum_{j=1}^i{i\choose j}(\pm 16p^h)^j\\
    &=&p\pm 16p^h n(n-1)/2+\sum_{i=1}^{n-1}\sum_{j=2}^i{i\choose j}(\pm 16p^h)^j
    \end{eqnarray*}

Since $n=|p|$, one has $m-h=1$. If we consider the positive sign, the above equality is $p=p+L$, and $L>0$, a contradiction. Otherwise, the above equality is
\begin{eqnarray*}
-2p&=&- 16p^h p(p-1)/2+\sum_{i=1}^{n-1}\sum_{j=2}^i{i\choose j}(- 16p^h)^j
\end{eqnarray*}
a contradiction.
\end{myitemize}
Assume $n=2$. Then $16p^m=y^2-1=(y-1)(y+1)$. If $p\mid\gcd(y-1,y+1)$, then $p=2$ (a contradiction as then we will have a Catalan's solution, $y^2-2^{m+4}=1$). Otherwise, we either have $y+1=\pm2^{\alpha}p^m$ and $y-1=\pm2^{4-\alpha}$, so
$y\in\{17,9,5,3,2,0,-1,-3,-7,-15\}$ with $p^m=5,3,14$ corresponding to $y=9,-7,-15$, or $y-1=\pm2^{\alpha}p^m$ and $y+1=\pm2^{4-\alpha}$, so $y=7, -17,-9$ with $p^m=3, 18,5$. These solutions are rejected.
\end{Proof}

\begin{Lemma}[Lemma 5.5, \cite{CaoChu}]
\label{lem:x2+16=yn}
The only positive integer solutions $(x,y,h,n)$ to the equation
\[ x^2+2^h= y^n,\;n>1,\;y\textrm{ odd},\;h>2\]
are  $(x,y,h,n)=(7,3,5,4)$ and $(x,y,n)=(2^{h-2}-1,2^{h-2}+1,2)$.
\end{Lemma}

Now we use some techniques in elementary number theory to find the integer solutions of some Diophantine equation.

\begin{Lemma}
\label{lem:x2-125=4yl}
The integer solutions $(x,y,l),\;l>1,y>0,$ to the Diophantine equation \begin{eqnarray}\label{eq3} x^2-125=\pm4y^l\end{eqnarray} are
\begin{eqnarray*}
\{(\pm15,5,2), (\pm63,31,2),(\pm11,1,l),(\pm5,5,2),(\pm25,5,3)\},
\end{eqnarray*}
or $l$ is odd, $5\nmid x$, and $x^2-125=4y^l$. 
\end{Lemma}
\begin{Proof}
We consider many possibilities:
\begin{myitemize}
\item[i.] $x^2-125=-4y^l$: By investigating perfect squares of the form $125-4\lambda$, one find the following possible solutions: \[\{(\pm11,1,l),(\pm9,11,1),(\pm7,19,1),(\pm5,5,2), (\pm3,29,1),(\pm1,31,1)\}\]
\item[ii.] $l=2k$ and $x^2-125=4y^l$: Then we can write $(x-2y^k)(x+2y^k)=125$. Therefore, we can assume that $(x-2y^k)\in\{\pm 1,\pm5,\pm25,\pm125\}$. Consequently, $l=2$ and we have

    $$\begin{array}{|c|c|c|c|c|c|c|c|}
    \hline
      x-2y^k &\pm1&\pm 5&\pm25&\pm125\\
      \hline
(x,y)& (\pm63,\pm31)& (\pm15,\pm5) & (\pm15,\mp5)&(\pm63,\mp31)\\
      \hline
    \end{array}$$

   \item[iii.] $l$ is odd, $x^2-125=4y^l$, and $5\mid x$: Then $5\mid y$. Since $l\ge 3$, one has $25\mid x$. In fact, $l=3$ and $5\mid\mid y$. Dividing by $125$, one has $5(x/25)^2-1=4(y/5)^3$. If $(x,y)$ is an integer solution to the latter equation, then $(X,Y)=\left(100(x/25),20(y/5)\right)$ is a solution to $X^2-k=Y^3,k=2000$. In \cite{MordellGebel}, all Mordell's equations with $|k|\le 10000$ were solved in $\Z$. In fact, the only solutions of $X^2-2000=Y^3$ are $(\pm100,20)$ and $(\pm44,-4)$. Therefore, the only integer solution of $x^2-125=4y^l,$ $l$ odd, and $5\mid x$ is $(\pm25,5,3)$.

    \item[iv.] $l$ is odd, $x^2-125=4y^l$, and $5\nmid x$: We notice that in (\ref{eq3}), both $x$ and $y$ are odd. In fact, this equation has a finite number of integer solutions. 
    \end{myitemize}
\end{Proof}
\begin{Corollary}
\label{cor:s2-11s-1=yl}
 The only integer solutions $(s,y,l)$ of the Diophantine equation \[s^2-11s-1=\pm y^l\] where $|s|,y>0$ are prime powers, $l>1$, are
 \begin{eqnarray*}\{(13,5,2), (-2,5,2), (37,31,2),  (11,1,l), (8, 5, 2), (3,5,2), (-7,5,3)\},\end{eqnarray*}
 or $l$ is odd and $s^2-11s-1=y^l$.
 \end{Corollary}
 \begin{Proof}
 After completing the square, one has
\begin{eqnarray*}
x^2-125=\pm4y^l,\textrm{ where }x=2s-11.
\end{eqnarray*}
According to Lemma \ref{lem:x2-125=4yl}, we obtain the above triples. We observe that the triples $(-26, 31,2)$ and $(18,5,3)$, corresponding to $(x,y,l)=(-63,31,2)$ and $(25,5,3)$ are rejected, because $|s|$ is not a prime power.
\end{Proof}

\section{Elliptic curves with rational $n$-torsion points}

Let $E/\Q$ be an elliptic curve with minimal discriminant $\Delta_E$ and conductor $N_E$. Given that $P\in E(\Q)[m],\;m\ge4$, there exist $b,c\in\Q$ such that the following Weierstrass equation defines an elliptic curve $E_{b,c}$ isomorphic to $E$
\begin{eqnarray}\label{equniversal}
E_{b,c}:y^2+(1-c)xy-by=x^3-bx^2
\end{eqnarray}
 with the image of $P$ being $(0,0)$. The invariants $c_4(b,c)$ and $\Delta(b,c)$ of $E_{b,c}$ are as follows:
\begin{eqnarray*}
c_4(b,c)&=&16b^2+8b(1-c)(c+2)+(1-c)^4\\
\Delta(b,c)&=&b^3\left(16b^2-b(8c^2+20c-1)-c(1-c)^3\right)
\end{eqnarray*}

By taking $m$ to be an integer in $\{4,5,6,7,8,9,10,12\}$, one obtains an explicit relation between $b,c$, see \S 2 of \cite{LorenziniTamagawa}.

\subsection{Case $n=4$} Assuming that $P\in E(\Q)[4]$, one has that $c=0$ in (\ref{equniversal}). We set $\lambda:=b$. The following Weierstrass equation describes $E$:
\[E:y^2+xy-\lambda y=x^3-\lambda x^2\]

Assume $\displaystyle\lambda=\frac{s}{t},\;s,t\in\Z,\;\gcd(s,t)=1.$ We obtain an integral Weierstrass equation describing $E$ using the following change of variables $x\mapsto x/t^2,\;y\mapsto y/t^3$. This integral equation is
\begin{eqnarray*}
E&:&y^2+t xy-st^2y=x^3-stx^2
\end{eqnarray*}
with the following invariants
\begin{eqnarray*}
\Delta_E&=& s^4t^7(16s+t)\\
c_4&=&t^2(16 s^2 + 16 s t + t^2)\\\;c_6&=&-t^3(-64 s^3 + 120s^2 t + 24 s t^2 + t^3)
\end{eqnarray*}

\begin{Theorem}
\label{thm:4}
Let $E/\Q$ be an elliptic curve such that $ E(\Q)[4]\ne\{0\}$. Assume moreover that $N_E=p q$ where $p\ne q$ are primes. Then $|\Delta_E|=p^{\alpha}q^{\beta}$ is given as follows:
\begin{eqnarray*}
  2^4\times 3,\; 2^4\times 5,\; 2^4\times 3^7,\;2^8\times7,\;2^8\times 3^2,\;2^8\times 7^7,\;3^2\times 7,\;3^2\times 5^2,\textrm{ and,}
\end{eqnarray*}
{\scriptsize$$\begin{array}{|c||c|c|c|c|c|c|c|c|c|}
    \hline
    |\Delta_E|&2^{2k+4}p&2^{2k+4} p^4&2^{4k}p&2^{4k}p^7&p^4q^b&p^4q^{7b}&p^{4k}q&p^{4k}q^7&p^{2k}q\\
    \hline
    p,q&\multicolumn{2}{|c|}{p=2^{k-4}\pm1,\;k\ge 4}&\multicolumn{2}{|c|}{p=2^{k+4}\pm1,\;k>0}&\multicolumn{2}{|c|}{16p\pm1=q^b,\;b>1}&\multicolumn{2}{|c|}{q=16p^k\pm1,\;k>1}&q=p^{2k}+16,\;k>0\\
    \hline
    \end{array}$$}
\end{Theorem}
\begin{Proof}
Let $s,t\in\Z$ be such that $E$ is given by the following Weierstrass equation \[E:y^2+t xy-st^2y=x^3-stx^2,\textrm{ where }\Delta_E=s^4t^7(16s+t)\]

One has $\gcd(s,t)=\gcd(s,16s+t)=1,$ and $\gcd(t,16s+t)=2^k,$ where $0\le k\le 4$, otherwise $2^{k-4}\mid s$, which is a contradiction. In fact, if $k>1$, then the Weierstrass equation is not minimal at $2$. Moreover, if $\ord_p(t)$ is odd, then $E$ has additive reduction at $p$.

We first treat the case that $|\Delta_E|=2^ap^b$ where $p$ is an odd prime, and $a,b>0$. Given $s$ and $t$, the following table gives the possible values for $\Delta_E=s^4t^7(16s+t)=2^{\alpha}p^{\beta}$. Observe that the table includes all possible values for $\Delta_E$, even when $E$ has additive reduction at some prime divisor of $t$.

{\footnotesize$$\begin{array}{|c|c|c|c|c|}
    \hline
    \textrm{\backslashbox{$|t|$}{$|s|$}}&2^m&2^mp^n&p^n&1\\
    \hline
    2^k,k>4&-&-&-2^5\times 3^8,2^3\times 3^4&2^5\times 3^2,2^3\times 3\\
    &&&\pm p^{4}2^{7k+4},p=2^{k-4}\pm1&\pm2^{7k+4} p,p=2^{k-4}\pm1\\
    \hline
    2^k,1\le k\le4&-&-&-&2^8\times 3^2,-2^8\times 7\\&&&&-2^4\times 3,2^4\times 5\\
    \hline
    2^kp^l,k>4,&-&-&-&-\\
    \hline
    2^kp^l,1\le k\le4&-&-&-&2^8\times 3^2,\pm2^8\times 7^7\\
    &&&&\pm2^4\times 3^7, 2^4\times 5^7\\
    \hline
    p^l&\pm2^{4m}p^7,\;p=2^{4+m}\pm1&-&-&-\\
    \hline
    1&\pm2^{4m}p,\;p=2^{4+m}\pm1&-&-&-    \\
      \hline
    \end{array}$$}
When $(|s|,|t|)=(2^m,p^l)$, we need to find solutions to $|16s+t|=|2^{4+m}\pm p^l|=1$. Therefore, we either have the unique Catalan solution or $l=1$, see Proposition \ref{prop:catalan}. The same reasoning and coprimality give the remaining possible values of $\Delta_E$ in the above table. Recall that if $\ord_p(t)>1$, then $E$ is not minimal at $p$, and we should consider $\ord_p(\Delta_E)$ mod $12$.

Now we assume $2\nmid N_E$. So without loss of generality we can assume that $\gcd(t,16s+t)=1$. Therefore, at least one of $|s|,|t|,|16s+t|$ is $1$.

\textbf{Case $|s|=1$:} Then $|\Delta_E|=|t^7(t\pm16)|$. Assuming $N_E=pq$, one observe that if $|t|=p^{a},a>1$, then $E$ is not minimal at $p$. In fact, if $a\equiv 1$ mod $2$, then after minimizing $E$ we obtain that $\ord_p(c_4)=2\ne0$ and hence $E$ has additive reduction at $p$ contradicting the fact that $p||N_E$. Therefore, we assume $|t|=p^{2a}$ and $|\Delta_E|=p^{2a}|p^{2a}\pm16|=p^{2a}q^b$. Thus we need to solve $|p^{2a}\pm16|=q^b$. Lemma \ref{lem:x2+16=yn} gives the solution $(p,a,q,b>1)=(\pm 3,1,\pm5,2)$ to the equation $p^{2a}+16=q^b$, and $\Delta_E=3^2\times 5^2$. A simple factorization argument shows that the only solution to $p^{2a}-16=q^b$ is $(\pm5,1,\pm3,2)$, and $\Delta_E=3^2\times 5^2$. The solution to $16-p^{2a}=q^b$ is $(\pm3,1,7,1)$, with $\Delta_E=3^2\times 7.$

\textbf{Case $|t|=1$:} Then $|\Delta_E|=|s^4(16s\pm1)|=p^{4a}q^b$, in other words, $|16s\pm1|=|16p^a\pm1|=q^b$. According to Lemma \ref{lem:16pm+1=ql}, $a=1$ or $b=1$, and so $|\Delta_E|=p^{4}q^b$ or $p^{4a}q$ respectively.

\textbf{Case $|16s+t|=1$:} Then we want to solve $16p^m=q^n\pm1$, and $\Delta_E=s^4t^7$. Again, according to Lemma \ref{lem:16pm+1=ql}, $|\Delta_E|=p^{4m}q^7$ or $p^4q^{7n}$.
\end{Proof}

\begin{Theorem}
Let $E/\Q$ be an elliptic curve such that $E(\Q)[4]\ne\{0\}$. Assume moreover that $N_E=pq$ where $p\ne q$ are primes. Then $|\Delta_E|<N_E^{32}$. In particular, $E/\Q$ satisfies Szpiro's conjecture.
\end{Theorem}
\begin{Proof}
We only need to check that $|\Delta_E|<N_E^{32}$ for the values of $|\Delta_E|$ given in Theorem \ref{thm:4}. This is straightforward for the first row of possible values of $|\Delta_E|$ appearing in Theorem \ref{thm:4}. In fact, $|\Delta_E|\le N_E^8$.

Now we are going to verify that $|\Delta_E|<N_E^{32}$ for the values of $|\Delta_E|$ in the table of Theorem \ref{thm:4}.
 \begin{eqnarray*}
 |\Delta_E|&=&2^{2k+4}p^{m}=2^{12}2^{2(k-4)}(2^{k-4}\pm1)^m<2^{12}[2\times(2^{k-4}\pm1)]^2(2^{k-4}\pm1)^m\\
 &<&2^{14}(2^{k-4}\pm1)^{2+m}<2^{10}(2^{k-4}\pm1)^{6+m}\le N_E^{10},\textrm{ where }m=1,4\\
 |\Delta_E|&=&2^{4k}p^m=2^{4k}(2^{k+4}\pm1)^m<(2^{k+4}\pm1)^4(2^{k+4}\pm1)^m=p^{m+4}<N_E^{11},\textrm{ where }m=1,7\\
 |\Delta_E|&=&p^{4}q^{mb}=p^4(16p\pm1)^{m}< p^4\times p^{4m}=p^{4m+4}<N_E^{4m+4},\textrm{ where }m=1,7\\
 |\Delta_E|&=&p^{4k}q^{m}=p^{4k}(16p^k\pm1)^m<(16p^k\pm1)^4(16p^k\pm1)^m=q^{m+4}<N_E^{m+4},\textrm{ where }m=1,7\\
 |\Delta_E|&=&p^{2k}q=p^{2k}(p^{2k}+16)<(p^{2k}+16)^2=q^2<N_E^2
 \end{eqnarray*}
\end{Proof}

\subsection{Case $n=5$} Assuming that $P\in E(\Q)[5]$, one has that $b=c$ in (\ref{equniversal}). Set $\lambda=b$. Then the following Weierstrass equation describes $E$:
\[E:y^2+(1-\lambda)xy-\lambda y=x^3-\lambda x^2\]

Assume $\displaystyle\lambda=\frac{s}{t},\;s,t\in\Z.$ We can obtain an integral Weierstrass equation describing $E$ using the following change of variables $x\mapsto x/t^2,\;y\mapsto y/t^3$. This integral equation is
\begin{eqnarray*}
E&:&y^2+(t-s)xy-st^2y=x^3-stx^2
\end{eqnarray*}
where the invariants of $E$ are given by
\begin{eqnarray*}
\Delta_E&=& s^5t^5(s^2-11st-t^2)\\
c_4&=&24st^2(-s + t) + (s^2 - 6 s t + t^2)^2
\end{eqnarray*}

\begin{Theorem}
\label{thm:5}
Let $E/\Q$ be an elliptic curve such that $E(\Q)[5]\ne\{0\}$. Assume moreover that $N_E=p^\alpha q^{\beta}$ where $p\ne q$ are primes, and $\alpha,\beta>0$. Then $|\Delta_E|=p^aq^b$ is given as follows:
\[2^5\times 5^2, \;2^{15}\times 5^2,\; 3^5\times 5^2,\; 13^5\times 5^2,\;37^5\times 31^2,\;7^5\times 5^3, \;p^{5k} q^{2l+1}\]
\end{Theorem}
\begin{Proof}
As we saw above there exist $s,t\in\Z$ such that $E$ is given by the equation
\[E:y^2+(t-s)xy-st^2y=x^3-stx^2\]
The assumption that $N_E$ is a product of two distinct prime powers together with the fact that $\gcd(t,s)=\gcd(s,s^2-11st-t^2)=\gcd(t,s^2-11st-t^2)=1$ imply that at least one of $|s|,|t|,|s^2-11st-t^2|$ is $1$.

\textbf{Case $|t|=1$:} Then one has $\Delta_E=s^5(s^2\mp11s-1)$. Consequently $s$ and $s^2\mp11s-1$ are both prime powers. Now we are going to solve the Diophantine equation $\displaystyle s^2\mp11s-1=\pm y^l$ and spot out the integer solutions $(s,y,l)$ where $s,y$ are prime powers. Completing the square, we need to find the integer solutions of \[x^2-125=\pm4y^l,\textrm { wehre }x=2s\mp11\] The solutions of the latter Diophantine equation is given in Lemma \ref{lem:x2-125=4yl}. In fact, we obtain the following table:

$$\begin{array}{|c|c||c|c|}
    (x,s,y,l)&\Delta_E&(x,s,y,l)&\Delta_E\\
    (\pm15,\pm13,5,2)&\pm13^5\times 5^2&(\pm25,\pm7,5,3)&\pm7^5\times 5^3\\
    (\pm15,\pm2,5,2)&\pm2^5\times 5^2&(\pm5,\pm2^3,5,2)&\pm2^{15}\times 5^2\\
    (\pm63,\pm37,31,2)&\pm37^5\times 31^2&(\pm5,\mp3,5,2) &\mp3^5\times 5^2\\
    (\pm63,\pm26,31,2)&\pm26^5\times 31^2&(x,s,q,1),|s|=p^k&\pm p^{5k}q\\
    \end{array}$$
or $|\Delta_E|=p^{5k}q^{2l+1}$ for finitely many possible primes $q$.  

\textbf{Case $|s|=1$:} Then $\Delta_E= t^5(1\mp11t-t^2)$. Similarly, we need to solve the Diophantine equation $125-x^2=\pm4y^l$ where $x=2t\pm11$. In fact, we obtain the same values given in the above table.

\textbf{Case $|s^2-11st-t^2|=1$:} Thus $|s|=p^m,|t|=q^n$ and $|\Delta_E|=p^{5m}q^{5n}$. Now we complete the square and have that $|(2p^m\mp 11q^n)^2-125q^{2n}|=4$. Any solution to the latter equation will yield a solution to the Diophantine equation $x^2-125y^{2}=\pm4$, where $x=2p^m\mp11q^n$ and $y=q^n$. We will start solving $x^2-125y^2=-4$ which is a Pell's equation for which we have the solution $(11,1)$. Thus any other solution $(x,y)$ is given by \[\frac{x+y\sqrt{125}}{2}=\pm\left(\frac{11+\sqrt{125}}{2}\right)^k,\]
see for example Proposition 6.3.16 in \cite{Cohentools}. One has
\begin{eqnarray*}
\frac{y}{2}&=&\frac{\pm1}{2^k}\left({k\choose 1}11^{k-1}+{k\choose 3}11^{k-3}\times 125+\ldots+{k\choose 1}11\times125^{(k-2)/2}\right),\textrm{ if }k\textrm{ is even}\\
\frac{x}{2}&=&\frac{\pm1}{2^k}\left(11^k+{k\choose2}11^{k-2}\times 125+{k\choose 4}11^{k-4}\times 125^2+\ldots+{k\choose 1}11\times 125^{(k-1)/2}\right),\textrm{ if } k \textrm{ is odd}
\end{eqnarray*}
Consequently, if $k$ is even, then $q=11$. Similarly, if $k$ is odd, then $11\mid x$ and $p=11$. If $q=11$, we have $|s^2\mp11^{n+1}s-11^{2n}|=1$, and $\displaystyle s=\frac{1}{2}\left(\pm11^{n+1}\pm\sqrt{11^{2n+2}-4(-11^{2n}\mp 1)}\right)$. This implies that there is a $\lambda\in\Z$ such that $\lambda^2=125\times 11^{2n}\pm4$. However, $\lambda^2+4=125\times 11^{2n}$ is not solvable by considering it mod $11$ as then $\lambda^2\equiv-4$ mod $11$ which contradicts that the Legendre symbol $\displaystyle \left(\frac{-4}{11}\right)=-1$. Moreover, $\lambda^2-4=125\times 11^{2n}$ can be shown to be non-solvable because $\lambda-2$ and $\lambda+2$ are coprime. Thus $\lambda-2\in\{\pm1,\pm125,\pm11^{2n},\pm125\times 11^{2n}\}$, and the only possible value for $|s|$ is $122$ which is not a prime power.
An identical argument holds when $p=11$.

Now we solve the equation $x^2-125y^2=4$. When $q$ is odd, one has that $x-2$ and $x+2$ are coprime. Following the same argument in the previous paragraph, $|s|$ cannot be a prime power. Now assume $q=2$. We divide by $4$ and obtain the new Diophantine equation $x^2-125\times 2^{2n-2}=1$, or $(x-1)(x+1)=125\times 2^{2n-2}$, and $x-1\in\{\pm2^{2n-h-2}\times 125,\,2^{2n-h-2}:h>0\}.$ Consequently, either $x+1=\pm2^h=\pm2^{2n-h-2}\times 125+2$, or $x+1=\pm2^h\times 125=2^{2n-h-2}+2$. Therefore, $h=1$ or $2n-h-1=1$ which is a contradiction in both cases.
\end{Proof}

\begin{Theorem}
Let $E/\Q$ be an elliptic curve such that $ E(\Q)[5]\ne\{0\}$. Assume moreover that $N_E=p^\alpha q^{\beta}$ where $p\ne q$ are primes, and $\alpha,\beta>0$. Then $|\Delta_E|<N_E^{6}$. In particular, $E/\Q$ satisfies Szpiro's conjecture.
\end{Theorem}
\begin{Proof}
We check that $|\Delta_E|<N_E^6$ for the possible values of $\Delta_E$ given in Theorem \ref{thm:5}. In fact, this is clear except when $|\Delta_E|=p^{5k}q$.

In the proof of Theorem \ref{thm:5}, we observe that $q=\mid p^{2k}\mp11p^k-1\mid$. In fact, $ p^{2k}\mp11p^k-1> p^k$ when $p^k>11$. In the latter case
\begin{eqnarray*}
N_E^6&=&p^6(p^{2k}\mp11p^k-1)^{6}= p^6 (p^{2k}\mp11p^k-1)^{5}(p^{2k}\mp11p^k-1)\\&>& p^{6+5k}(p^{2k}\mp11p^k-1)> p^{5k}q=\mid\Delta_E\mid
\end{eqnarray*}
  We are left with treating a finite number of possible cases, namely
    \[p^k\in\{2,3,4,5,7,8,9,11\}\]
  Straight forward calculations show that for all these cases if $N_E$ is a product of two distinct prime powers, then $\mid\Delta_E\mid\le N_E^6$.
\end{Proof}

\subsection{Case $n=6$}
Let $P\in E(\Q)[6]$. There exists a $\lambda\in\Q$ such that the following Weierstrass equation describes $E$: \[y^2+(1-\lambda)xy-\lambda(\lambda+1)y=x^3-\lambda(\lambda+1)x^2\]
Assuming $\lambda=s/t$, $\gcd(s,t)=1$, we may use the transformation $\displaystyle x\mapsto x/t^2,\;y\mapsto y/t^3$ to obtain the following Weierstrass equation
\begin{eqnarray*}
y^2+(t-s)xy-(t^2s+ts^2)y&=&x^3-(st+s^2)x^2
\end{eqnarray*}
where the invariant of $E$ are given by
\begin{eqnarray*}
\Delta_E&=&s^6t^2(s+t)^3(9s+t)\\
c_4&=&(3s+t)(3s^3+3s^2t+9st^2+t^3)
\end{eqnarray*}

\begin{Theorem}
\label{thm:6}
Let $E/\Q$ be an elliptic curve such that $ E(\Q)[6]\ne\{0\}$. Assume moreover that $N_E=p^{\alpha} q^{\beta}$ where $p\ne q$ are primes, and $\alpha,\beta>0$. Then $\Delta_E$ is given as follows:
\begin{eqnarray*}
2\times 7^2,\; -2^2\times 7,\;2^3\times7^6,\;2^4\times 5,\;-2^4\times 3^3,\;2^6\times 17,\;-2^6\times 7^3,\;2^8\times 3^3,\;-2^8\times5^2
\end{eqnarray*}

In particular, $|\Delta_E|<N_E^6$, and Szpiro's conjecture holds for $E$.
\end{Theorem}
\begin{Proof}
Let $s,t\in\Z$, $\gcd(s,t)=1$, be such that the following Weierstrass equation describes $E$
\begin{eqnarray*}
y^2+(t-s)xy-(t^2s+ts^2)y&=&x^3-(st+s^2)x^2
\end{eqnarray*}
where \[\gcd(s,t)=\gcd(s,s+t)=\gcd(t,s+t)=\gcd(s,9s+t)=1\]
and \[\gcd(t,9s+t)\mid 9,\;\gcd(s+t,9s+t)\mid 8\]

\textbf{Case i.} Assume $\gcd(t,9s+t)=\gcd(s+t,9s+t)=1$. Then at least two of $|s|,|t|,|s+t|,|9s+t|=1$. If $s=t=\pm1,$ then $s+t=\pm2$, $9s+t=\pm10$, and $\Delta_E=2^4\times 5$. When $s=\pm1,s+t=\mp1$, one has $t=\mp 2$, $9s+t=\pm7$, and $\Delta_E=-2^2\times7.$ When $|s|=|9s+t|=1$, one has that $(s,t,9s+t,s+t)\in\{(\pm1,\mp8,\pm1,\mp7),(\pm1,\mp10,\mp1,\mp9)\}$, and the first quadruple yields $\Delta_E=-2^6\times7^3.$ When $t=\pm 1, s+t=\mp1$, one has $s=\mp2$, $9s+t=\mp17$, and $\Delta_E= 2^6\times 17$.
The possibilities $|t|=|9s+t|=1$ and $|s+t|=|9s+t|=1$ are rejected.

\textbf{Case ii.} Assume $\gcd(t,9s+t)=3^k$ and $\gcd(s+t,9s+t)=2^l$ where $1\le k\le2,\le l\le3$. Then $|t|=3^f,\;|s+t|=2^g$ and $|s|=1$. In other words, $|1\pm3^f|=2^g$. According to Proposition \ref{prop:catalan}, the only integer solutions of the latter equation will yield the following set of quadruples \[(s,t,s+t,9s+t)\in\{(\pm1,\pm3,\pm2^2,\pm12),(\pm1,\mp3,\mp2,\pm6),(\pm1,\mp3^2,\mp2^3,0)\}.\] The first quadruple gives $\Delta_E=2^8\times 3^3$, whereas the second gives $\Delta_E=-2^4\times3^3$.

\textbf{Case iii.} Now assume $\gcd(s+t,9s+t)=2^l,\;1\le l\le3$, and $\gcd(t,9s+t)=1$. The following table provides the possible values of the discriminant of $E$:

{\footnotesize$$\begin{array}{|c|c|c|}
    \hline
    \textrm{\backslashbox{$|s+t|$}{$|9s+t|$}}&2^f&2^fp^m\\
    \hline
    2^g&-2^8\times5^2,\;2^{13}\times 7^2,\;2^{15}\times7^6&-\\
    \hline
    2^gq^n&2^4\times 5&-    \\
      \hline
    \end{array}$$}
If $(|s+t|,|9s+t|)=(2^f,2^g)$, then $2^g=|9s+t|=|8s\pm2^f|$. It follows that either $f=g=2$ or $\min(f,g)=3$. We obtain the following quadruples $(s,t,s+t,9s+t)\in\{(\pm1,\pm7,\pm2^3,\pm2^4),(\pm1,\mp5,\mp2^2,\pm2^2),(\pm7,\pm1,\pm2^3,\pm2^6)\}$. The first and third quadruples give non-minimal elliptic curves and so we have to minimize them. If $(|s+t|,|9s+t|)\in\{(2^fp^m,2^g),(2^f,2^gq^n)\}$, then $|s|=|t|=1$ and the second pair gives $s+t=\pm2,9s+t=\pm10$.

\textbf{Case iv. } Assume $\gcd(t,9s+t)=3^l,\;l\in\{1,2\}$, and $\gcd(s+t,9s+t)=1.$ In fact, the only prime divisor of $t$ and $9s+t$ is $3$, since otherwise $|s|=|s+t|=1$ where $t$ is divisible by $3$. Therefore, $(|t|,|9s+t|)=(3^f,3^g)$ where $\min(f,g)\in\{1,2\}$, and $|\Delta_E|=s^63^{2f}(s+t)^3(\pm3^g)$. Either $|s|$ or $|s+t|$ is $1$. If $|s|=1$, then $3^g=|9s+t|=|\pm9+3^f|$ and there is no integer $s$ satisfying the latter equalities. If $|s+t|=1$, then $3^g=|9s+t|=|9(s+t)-8t|=|9\pm8\times3^f|$. The only quadruple $(s,t,s+t,9s+t)$ satisfying the latter equalities under the condition that $\min(f,g)\in\{1,2\}$ is $(\pm10,\mp3^2,\pm1,\pm3^4)$, but then $\Delta_E$ has three distinct prime divisors.
\end{Proof}

\subsection{Case $n=7$}
Let $P\in E(\Q)[7]$. Then there exists a $\lambda\in\Q$ such that the following Weierstrass equation describes $E$:
\[y^2+(1-\lambda(\lambda-1))xy-\lambda^2(\lambda-1)y=x^3-\lambda^2(\lambda-1)x^2\]
\begin{Theorem}
\label{thm:7}
Let $E/\Q$ be an elliptic curve such that $ E(\Q)[7]\ne\{0\}$. Assume moreover that $N_E=p^{\alpha} q^{\beta}$ where $p, q$ are distinct primes, and $\alpha,\beta>0$. Then $\Delta_E=-2^7\times 13$.
In particular, $|\Delta_E|<N_E^3$, and Szpiro's conjecture holds for $E$.
\end{Theorem}
\begin{Proof}
Assuming $\lambda=s/t$, $\gcd(s,t)=1$, we use the transformation $\displaystyle x\mapsto x/t^4,\;y\mapsto y/t^6$ to obtain the following integral Weierstrass equation describing $E$
\begin{eqnarray*}
y^2+(t^2-s^2+st)xy-s^2(st^3-t^4)y=x^3-(s^3t-s^2t^2)x^2\\
\end{eqnarray*}
with invariants
\begin{eqnarray*}
\Delta_E&=&s^7t^7(s-t)^7(s^3-8s^2t+5st^2+t^3)\\
c_4&=&(s^2-st+t^2)(s^6-11s^5t+30s^4t^2-15s^3t^3-10s^2t^4+5st^5+t^6)
\end{eqnarray*}
We set $k=s^3-8s^2t+5st^2+t^3$. Then
\begin{eqnarray*}\gcd(s,t)=\gcd(s,s-t)=\gcd(s,k)=\gcd(t,s-t)
=\gcd(t,k)=\gcd(s-t,k)=1.
\end{eqnarray*}
 Since $\Delta_E$ has only two distinct prime divisors, then at least two of $|s|,|t|,|s-t|,|k|$ are ones. Indeed, if two of $|s|,|t|,|s-t|$ are ones, then the only corresponding discriminant is $\Delta_E=-2^7\times 13.$

If $|s|=|k|=1$, then $|\pm1-8t\pm 5t^2+t^3|=1$ and either $t=0$, or $t=s=\pm1$ which yields $s-t=0$. The same holds if $|t|=|k|=1$ or $|s-t|=|k|=1$.
\end{Proof}

\subsection{Case $n=8$}
Let $P\in E(\Q)[8]$. Then there exists a $\lambda\in\Q$ such that $E$ is described by the following Weierstrass equation:
\[y^2+\left(1-\frac{(2\lambda-1)(\lambda-1)}{\lambda}\right)xy-(2\lambda-1)(\lambda-1)y=x^3-(2\lambda-1)(\lambda-1)x^2\]
\begin{Theorem}
\label{thm:8}
Let $E/\Q$ be an elliptic curve such that $E(\Q)[8]\ne\{0\}$. Assume moreover that $N_E=p^{\alpha} q^{\beta}$ where $p\ne q$ are primes, and $\alpha,\beta>0$. Then $\Delta_E=-2^{11}\times3^8$.
In particular, $|\Delta_E|<N_E^7$, and Szpiro's conjecture holds for $E$.
\end{Theorem}
\begin{Proof}
In the above Weierstrass equation we take $\lambda=s/t$ where $\gcd(s,t)=1$. Then we apply the transformation
\[x\mapsto\frac{x}{s^2t^2},\;y\mapsto\frac{y}{s^3t^3}\] to obtain the following integral Weierstrass equation
\[y^2-(t^2-4st+2s^2)xy-ts^3(s-t)(2s-t)y=x^3-s^2(s-t)(2s-t)x^2\] where
\begin{eqnarray*}
\Delta_E&=&s^8 t^2(s-t)^8(2s-t)^4(8s^2-8st+t^2)\\
c_4&=&16s^8 - 64 s^7 t + 224 s^6 t^2 - 448 s^5 t^3 +
    480s^4 t^4 - 288 s^3 t^5 + 96 s^2 t^6 - 16 s t^7 + t^8
\end{eqnarray*}
If $t$ is even, one has that $|s|$, $|t|$, $|s-t|$ and $|s-t/2|$ are pairwise coprime, so at least two of these are ones.
If $|s|=|s-t|=1$, then $s=\pm 1,t=\pm 2$ and $2s-t=0$. If $|s|=|s-t/2|=1$, then $s=\pm1,\;t=\pm4$, $s-t=\mp3$,\; $8s^2-8st+t^2=-8$ and $\Delta_E=-2^{11}\times 3^8$. If $s-t=\pm1,s-t/2=\mp1$, then $s=\mp3,t=\mp4,2s-t=\mp2,8s^2-8st+t^2=-8$ and $\Delta_E=-2^{11}\times3^8$.

If $t$ is odd, then at least two of $|s|$, $|t|$, $|s-t|$ and $|2s-t|$ are ones. The following table contains $|s t(s-t)(2s-t)(8s^2-8st+t^2)|$.

{\footnotesize$$\begin{array}{|c|c|c|c|c|}
    \hline
 {}&|s|=1&|t|=1&|2s-t|=1,\,t\textrm{ odd}\\
    \hline
    |s-t|=1&t\textrm{ even}&2\times3\times17&2\times3\times7\\
      \hline
      |2s-t|=1,\,t\textrm{ odd}&2\times 3\times7&0\\
      \cline{1-3}
      |t|=1&2\times 3\times 17\\
      \cline{1-2}
    \end{array}$$}

Therefore, the only elliptic curve with an $8$-torsion point and whose conductor is a product of two distinct prime powers is the one with $\Delta_E=-2^{11}\times 3^8$. We observe that its invariant $c_4=-2^4\times47$, so it has additive reduction over $\mathbb{F}_2$ and $N_E=2^2\times 3.$
\end{Proof}

\subsection{Case $n=9$}
Let $P\in E(\Q)[9]$. There exists a $\lambda\in\Q$ such that $E$ is described by the following Weierstrass equation:
\[y^2+\left(1-\lambda^2(\lambda-1)\right)xy-\lambda^2(\lambda-1)(\lambda^2-\lambda+1)y=x^3-\lambda^2(\lambda-1)(\lambda^2-\lambda+1)x^2\]
\begin{Theorem}
\label{thm:9}
Let $E/\Q$ be an elliptic curve such that $E(\Q)[9]\ne\{0\}$. Assume moreover that $N_E=p^{\alpha} q^{\beta}$ where $p\ne q$ are primes, and $\alpha,\beta>0$. Then $\Delta_E=-2^9\times 3^5$.
In particular, $|\Delta_E|<N_E^5$, and Szpiro's conjecture holds for $E$.
\end{Theorem}
\begin{Proof}
In the above Weierstrass equation we take $\lambda=s/t$ where $\gcd(s,t)=1$. Then we apply the transformation
\[x\mapsto\frac{x}{t^6},\;y\mapsto\frac{y}{t^9}\] to obtain the following integral Weierstrass equation describing $E$
\[y^2+\left(t^3-s^2(s-t)\right)xy-t^4s^2(s-t)(s^2-st+t^2)y=x^3-ts^2(s-t)(s^2-st+t^2)x^2\] where
\begin{eqnarray*}
\Delta_E&=&s^9t^9(s-t)^9(s^2-st+t^2)^3(s^3-6s^2t+3st^2+t^3)\\
c_4&=&(s^3-3s^2t+t^3)(s^9 - 9 s^8 t + 27 s^7 t^2 - 48 s^6 t^3 + 54 s^5 t^4 -
    45 s^4 t^5 + 27 s^3 t^6 - 9 s^2 t^7 + t^9)
\end{eqnarray*}
since
$s^2-st+t^2=(s-t)^2+st$, the first four factors $s,t,s-t,s^2-st+t^2$ of $\Delta_E$ are pairwise coprime. Therefore, at least two of the absolute values of these factors are ones. If $s=\pm1,t=\mp1$, then $s-t=\pm2$, $s^2-st+t^2=3$, $s^3-6s^2t+3st^2+t^3=\pm9$, $\Delta_E=-2^9\times 3^5$. If $|s|=|s-t|=1$, then $s=\pm1,t=\pm2,s-t=\mp1$, $s^2-st+t^2=3$, $s^3-6s^2t+3st^2+t^3=\pm9$, $\Delta_E=-2^9\times 3^5$. If $s=\pm1$, $-1=s^2-st+t^2=1\mp t+t^2$, then there is no integer $t$ satisfying the latter equalities. If $s=\pm1$, $|s^3-6s^2t+3st^2+t^3|=1$, then $t=\pm1.$ We will have the same results if we replace $|s|=1$ by $|t|=1$. If $|s-t|=1$ and $|s^2-st+t^2|=|(s-t)^2+st|=|1+st|=1$, then $st=-2$, a contradiction.

The only elliptic curve with a $9$-torsion point and whose conductor is a product of two distinct prime powers is the one with $\Delta_E=-2^{9}\times 3^5$. We observe that $\ord_3(c_4)>0$, so it has additive reduction over $\mathbb{F}_3$ and $N_E=2\times 3^2.$
\end{Proof}

\subsection{Case $n=10$}
Let $P\in E(\Q)[10]$. There exists a $\lambda\in\Q$ such that the following Weierstrass equation describes $E$:
\[y^2+\left(1+\frac{\lambda(\lambda-1)(2\lambda-1)}{(\lambda^2-3\lambda+1)}\right)\,xy-\frac{\lambda^3(\lambda-1)(2\lambda-1)}{(\lambda^2-3\lambda+1)^2}\,y=
x^3-\frac{\lambda^3(\lambda-1)(2\lambda-1)}{(\lambda^2-3\lambda+1)^2}\, x^2\]
\begin{Theorem}
\label{thm:10}
There exists no elliptic curve $E/\Q$ with $E(\Q)[10]\ne\{0\}$ and $N_E=p^{\alpha} q^{\beta}$ where $p\ne q$ are primes, and $\alpha,\beta>0$.
\end{Theorem}
\begin{Proof}
In the above Weierstrass equation we take $\lambda=s/t$ where $\gcd(s,t)=1$. Then we apply the transformation
\[x\mapsto\frac{x}{t^2(s^2-3st+t^2)^2},\;y\mapsto\frac{y}{t^3(s^2-3st+t^2)^3}\] to obtain the following integral Weierstrass equation
\begin{eqnarray*}y^2+\left[t(s^2-3st+t^2)+s(s-t)(2s-t)\right]xy-t^2s^3(s-t)(2s-t)(s^2-3st+t^2)y\\=x^3-ts^3(s-t)(2s-t)x^2\end{eqnarray*} where
\begin{eqnarray*}
\Delta_E=s^{10}t^5(s - t)^{10}(2s - t)^5 (4 s^2 - 2 s t - t^2)(s^2 - 3 s t + t^2)^2\\
\end{eqnarray*}
The factors $s,t,s-t,s^2-3st+t^2$ are pairwise coprime. Therefore, at least two of $|s|,|t|,|s-t|,|s^2-3st+t^2|$ are ones. The following table contains the product $|st(s - t)(2s - t)(4 s^2 - 2 s t - t^2)(s^2 - 3 s t + t^2)|$ when two of $|s|,|t|,|s-t|,|s^2-3st+t^2|$ are ones.

{\footnotesize$$\begin{array}{|c|c|c|c|c|c|}
    \hline
 {}&|s|=1&|t|=1&|s-t|=1&|2s-t|=1&|4s^2-2st-t^2|=1\\
    \hline
    |s^2-3st+t^2|=1&66& 870,66&0,66&870,66,0&0\\
      \hline
|4s^2-2st-t^2|=1&60&0&60&0\\
      \cline{1-5}
      |2s-t|=1&66&0&30\\
      \cline{1-4}
      |s-t|=1&0&66\\
      \cline{1-3}
      |t|=1&30\\
      \cline{1-2}
    \end{array}$$}

It is clear that $E$ cannot have a discriminant with only two distinct prime divisors.
\end{Proof}

\subsection{Case $n=12$}
Let $P\in E(\Q)[12]$. There exists a $\lambda\in\Q$ such that the following Weierstrass equation describes $E$:
\begin{eqnarray*}y^2+\left(1+\frac{\lambda(2\lambda-1)(3\lambda^2-3\lambda+1)}{(\lambda-1)^3}\right)\,xy-\frac{\lambda(2\lambda-1)(3\lambda^2-3\lambda+1)
(2\lambda^2-2\lambda+1)}{(\lambda-1)^4}\,y\\
=x^3-\frac{\lambda(2\lambda-1)(3\lambda^2-3\lambda+1)(2\lambda^2-2\lambda+1)}{(\lambda-1)^4}\, x^2\end{eqnarray*}
\begin{Theorem}
\label{thm:12}
There exists no elliptic curve $E/\Q$ with $E(\Q)[12]\ne\{0\}$ and $N_E=p^{\alpha} q^{\beta}$ where $p\ne q$ are primes, and $\alpha,\beta>0$.
\end{Theorem}
\begin{Proof}
In the above Weierstrass equation we take $\lambda=s/t$ where $\gcd(s,t)=1$. Then we apply the transformation
\[x\mapsto\frac{x}{t^2(s-t)^6},\;y\mapsto\frac{y}{t^3(s-t)^9}\] to obtain the following integral Weierstrass equation
\begin{eqnarray*}y^2+\left[t(s-t)^3+s(2s-t)(3s^2-3st+t^2)\right]xy-ts(s-t)^5(2s-t)(3s^2-3st+t^2)(2s^2-2st+t^2)y\\=x^3-s(s-t)^2(2s-t)(3s^2-3st+t^2)(2s^2-2st+t^2)x^2\end{eqnarray*} where
\begin{eqnarray*}
\Delta_E=s^{12}t^2(s - t)^{12}(2 s - t)^6(3 s^2 - 3 st + t^2)^4(2 s^2 - 2 s t + t^2)^3(6s^2- 6s t + t^2)\\
\end{eqnarray*}

If $t$ is even, then two of $|s|,|t|,|s-t|$ and $|s-t/2|$ are ones. If $t$ is odd, then two of $|s|,|t|,|s-t|$ and $|2s-t|$ are ones. In both cases, one finds that the product corresponding to $\Delta_E$ is either $0$ or has more than two prime divisors.
\end{Proof}

\hskip-18pt\emph{\bf{Acknowledgements.}}
I would like to thank Michael Bennett for spotting out a mistake in the proof of Theorem \ref{thm:5}. 
I would also like to thank Andrzej Dabrowski for sharing some of the results of his current work on elliptic curves with small conductor.

\bibliographystyle{plain}
\footnotesize
\bibliography{conductor}

\begin{thebibliography}{10}

\bibitem{BugeaudRoy}
Y.~Bugeaud, M.~Mignotte, and Y.~Roy.
\newblock On the diophantine equation $\displaystyle\frac{x^n-1}{x-1}=y^q$.
\newblock {\em Pacific {J}ournal of {M}athematics}, 193(2):257--268, 2000.

\bibitem{CaoChu}
Zhenfu Cao, Chuan~I Chu, and Wai~Chee Shiu.
\newblock The exponential {D}iophantine equation $ax^2+by^2=\lambda k^z$.
\newblock {\em TAIWANESE JOURNAL OF MATHEMATICS}, 12(5):1015--1034, August
  2008.

\bibitem{Cohentools}
Henri Cohen.
\newblock {\em Number {T}heory {V}olume {I}: {T}ools and {D}iophantine
  {E}quations}.
\newblock Graduate {T}exts in {M}athematics 239. Springer {V}erlag, 2007.

\bibitem{MordellGebel}
J.~Gebel, A.~Peth\H{o}, and H.~Zimmer.
\newblock On {M}ordell's equation.
\newblock {\em Compositio Mathematica}, 110(3):335--367, 2011.

\bibitem{Hadanoremarks}
T.~Hadano.
\newblock Remarks on the conductor of an elliptic curve.
\newblock {\em Proc. {J}apan {A}cad.}, 48:166--167, 1972.

\bibitem{Hadano}
T.~Hadano.
\newblock On the conductor of an elliptic curve with a rational point of order
  $2$.
\newblock {\em Nagoya {M}ath. {J}.}, 53:199--210, 1974.

\bibitem{LorenziniTamagawa}
D.~Lorenzini.
\newblock Torsion and {T}amagawa numbers.
\newblock {\em accepted for publication in {A}nn. {I}nst. {F}ourier}.

\bibitem{Ogg2power}
A.~P. Ogg.
\newblock Abelian curves of 2-power conductor.
\newblock {\em Proc. {C}amb. {P}hil. {S}oc.}, 62:143--148, 1966.

\bibitem{Oggsmallconductor}
A.~P. Ogg.
\newblock Abelian curves of small conductor.
\newblock {\em J. reine angew. Math.}, 226:205--215, 1967.

\bibitem{Alice}
A.~Silverberg.
\newblock Open questions in arithmetic algebraic geometry.
\newblock in Arithmetic Algebraic Geometry (Park City, UT, 1999), Institute for
  Advanced Study/Park City Mathematics Series 9, American Mathematical Society,
  Providence, RI, 2001.

\bibitem{sil1}
J.~Silverman.
\newblock {\em The arithmetic of elliptic curves}.
\newblock GTM 106. Springer-Verlag, New York, 1986.

\bibitem{Sil2}
J.~Silverman.
\newblock {\em Advanced topics in the arithmetic of elliptic curves}.
\newblock GTM 151. Springer-Verlag, 1995.

\end{thebibliography}
\end{document}